\documentclass[12pt]{article}
\usepackage[latin1]{inputenc}
\usepackage{amssymb,amsmath,amsfonts}
\newtheorem{theorem}{Theorem}

\newtheorem{lemma}[theorem]{Lemma}

\newcommand{\C}{\mathbb{C}}

\newcommand{\R}{\mathbb{R}}

\newcommand{\spa}{\mbox{span}}

\newcommand{\po}{{\hspace*{-1ex}}{\bf .  }}
\newcommand{\nab}{\tilde\nabla}

\renewcommand{\o}{\omega}
\def\bea{\begin{eqnarray*} }
\def\eea{\end{eqnarray*} }
\def\<{\langle}
\def\>{\rangle}
\def\n{\nabla}

\def\a{\alpha}

\def\e{\epsilon}
\def\span{\operatorname{span}}
\def\be{\begin{equation} }
\def\ee{\end{equation} }

\def\nab{\tilde\nabla}
\def\nap{\nabla^\perp}
\def\proof{\noindent{\it Proof:  }}
\def\qed{\ifhmode\unskip\nobreak\fi\ifmmode\ifinner\partial_z\alpha_2
\else\hskip5 pt \fi\fi\hbox{\hskip5 pt \vrule width4 pt
height6 pt  depth1.5 pt \hskip 1pt }}

\begin{document}

\title{A new class of austere submanifolds}
\author{M.\ Dajczer and Th.\ Vlachos}
\date{}
\maketitle

\begin{abstract}
Austere submanifolds of Euclidean space were introduced in 1982 by Harvey
and Lawson in their foundational work on calibrated geometries. In general, 
the austerity condition is much stronger than minimality since it express 
that the nonzero eigenvalues of the shape operator of the submanifold 
appear in opposite pairs for any normal vector at any point. Thereafter,
the challenging task of finding non-trivial explicit examples, other than
minimal immersions of Kaehler manifolds, only turned out submanifolds 
of rank two, and these are of limited interest in the sense that in this 
special situation austerity is equivalent to minimality. In this 
paper, we present the first explicitly given family of austere non-Kaehler 
submanifolds of higher rank, and these are produced from holomorphic data 
by means of a Weierstrass type parametrization.
\end{abstract}

After the celebrated paper by Harvey and Lawson \cite{HL} on calibrated 
geometries the classification of austere Euclidean submanifolds 
became a rather challenging task in submanifold theory. An isometric 
immersion $f\colon M^n\to\R^N$ of a Riemannian manifold $M^n$, $n\geq 2$, 
into Euclidean space is called \emph{austere} if the nonzero eigenvalues 
of the shape operator for any normal vector at any point appear in opposite 
pairs, or equivalently, if all odd degree symmetric polynomials on these 
eigenvalues vanish.  

The notion of austerity was introduced by Harvey and Lawson \cite{HL} in 
connection with the class of special Lagrangian submanifolds in complex 
Euclidean space $\C^N$ that are not only minimal but absolutely area 
minimizing. Given an isometric immersion $f\colon M^n\to\R^N$,  the 
embedding of its normal bundle $\psi\colon N_fM\to\R^N\oplus\R^N$ 
defined by
$$
\psi(\xi(x))=(f(x),\xi(x))
$$
is a Lagrangian submanifold of $\C^N\equiv\R^N\oplus\R^N$ with 
respect to the complex structure $J(X,Y)=(-Y,X)$. Then $\psi$ is 
special Lagrangian if and only if $f$ is austere.  

In the special case of a  submanifold $M^n$ in $\R^N$ of rank $\rho=2$, 
that is, when the kernel of the second fundamental form (called the 
\emph{relative nullity} subspace) of the submanifold has constant dimension 
(called the \emph{index of relative nullity}) $n-2$, we have that austerity 
and minimality are equivalent. Notice 
that $\rho$ is the rank of the Gauss map with values in the Grassmannian 
$G_{n,N}$ of oriented subspaces.  But for submanifolds for higher rank, 
the austerity condition is much more demanding than minimality. This 
makes it rather hard to find examples of austere submanifolds other 
than the obvious examples of holomorphic isometric immersions of Kaehler 
manifolds into $\C^N$. In fact, we know from \cite{DR} that for an 
isometric immersion of a Kaehler manifold into $\R^N$ to be austere 
it suffices to be minimal, but these immersions are always the ``real part" 
of a holomorphic one in $\C^N$.  

The quest to construct new examples of austere submanifolds was initiated 
by Bryant \cite{Br} who classified the rank two submanifolds of dimension 
three as well as a quite simple family of examples of higher dimension 
called generalized helicoids. Bryant showed that the interesting examples 
of dimension three are ``twisted cones" over minimal surfaces in spheres. 
As for dimension four, he provided a careful full pointwise description 
of the structures of all possible second fundamental forms. In a somehow 
dual parametric  form, Bryant's  construction in the three dimensional case 
was extended by Dajczer and Florit \cite{DF} to submanifolds of rank two of 
any dimension. Roughly speaking, they showed that these submanifolds are 
subbundles of the normal bundles of a class of Euclidean or spherical 
surfaces called elliptic that, in addition,  satisfy that the ellipses 
of curvature of a certain order are circles. But outside special cases, 
it is not known how to generate these surfaces. Finally, the four 
dimensional case was intensively studied by Ionel and Ivey \cite{II1}, 
\cite{II2} building on Bryant's algebraic results.  In particular, 
they obtained a non-parametric classification in the special case of 
the submanifolds ruled by planes.

In this paper, we take advantage of our results in \cite{DV} in order 
to characterize in an explicit parametric form a class of austere 
submanifolds $M^n$ in $\R^{n+2}$ of dimension $n\geq 4$ and rank $\rho=4$. 
Besides being the first non-trivial known examples, other than minimal 
Kaehler submanifolds, having any possible dimension and rank $\rho> 2$, 
what makes this new class of particular interest is that they are given 
in terms of a Weierstrass type parametrization depending on $n$ holomorphic 
functions on a domain. Consequently, the same is true for the special 
Lagrangian submanifolds that can be constructed from them as shown above.
\vspace{1ex}

Before stating our results, we first briefly recall some facts that
can be seen exposed with many details in \cite{DV}. In fact,  in the 
sequel we will make systematic use of results in that paper, sometimes 
without further referrence.
\vspace{1ex}

A substantial minimal surface $g\colon L^2\to\R^N$ is called 
\emph{$m$-isotropic}, $m\geq 1$, if at any point of $L^2$ all 
ellipses of curvature (defined below) until order $m$ are circles. 
Being substantial means that the surface is not contained in any 
proper affine subspace of $\R^N$, in fact, not even locally since
$g$ is real analytic.
It is well-known that $g\colon L^2\to\C^{N/2}\cong\R^N$ for $N$  
even is a holomorphic curve if and only if the ellipses of 
curvature of any order at any point are circles; for instance 
see \cite{CCh}.  

Any simply connected $m$-isotropic surface admits a Weierstrass type 
representation given in \cite{DG} based on results in \cite{CCh}.
In particular, any simply connected $2$-isotropic surface is obtained 
as follows: Start with a nonzero holomorphic map 
$\alpha_0\colon U\to\C^{N-4}$ on a domain $U\subset\C$ and define 
$\alpha_1\colon U\to\C^{N-2}$ by 
$$
\alpha_1=\beta_1\left(1-\phi_0^2,i(1+\phi_0^2),2\phi_0\right)
$$
where $\phi_0=\int_U^z\alpha_0dz$ and $\beta_1\neq 0$ is any holomorphic 
function. Define $\alpha\colon U\to\C^N$ 
by
$$
\alpha=\beta_2\left(1-\phi_1^2,i(1+\phi_1^2),2\phi_1\right)
$$
where $\phi_1=\int_U^z\alpha_1dz$ and $\beta_2\neq 0$ is any holomorphic 
function. If $\phi=\int_U^z\alpha dz$ then $g=\mbox{Re}\,\phi$ 
is a $2$-isotropic surface in $\R^N$.
\vspace{1ex}

It is easy to see that the above procedure yields examples of $2$-isotropic 
surfaces with complete metrics. For instance, see the construction at the 
final part of \cite{DG}.
\medskip

Let $g\colon L^2\to\R^{n+2}$ be a  $1$-isotropic oriented surface. Then 
let $\Lambda_g\subset N_gL$ be the vector subbundle of the normal 
bundle of $g$ with $(n-2)$-dimensional fibers
$$
\Lambda_g(u,v)=\left(\span\{g_u,g_v,g_{uu},g_{uv}\}\right)^\perp
$$
where $g=g(u,v)$ is parametrized in local isothermal coordinates. If 
$g=\mbox{Re}\,\phi$ is as above, then
$$
\Lambda_g=\alpha\wedge\alpha_z.
$$

It was shown in 
\cite{DV} that the dimension of $\Lambda_g(u,v)$ may fail to be $n-2$ 
only at isolated points and that the vector bundle extends smoothly to 
these points. Hence, from now on $\pi\colon\Lambda_g\to L^2$ denotes 
the extended vector bundle.

Let $F_g\colon\Lambda_g\to\R^{n+2}$ be the immersion associated to $g$
defined on $\pi\colon\Lambda_g\to L^2$ by
\be\label{param}
F_g(p,\xi)=g(p)+\xi,\;\;p=\pi(\xi).
\ee
In the sequel, we denote by $M^n$ the manifold $\Lambda_g$ when  endowed 
with the metric induced by $F_g$ and by $j\colon L^2\to M^n$ the immersion
in $M^n$ of the zero-section of $\Lambda_g$. We have by construction that 
$F_g$ is a $(n-2)$-ruled submanifold, and it is easily seen that $j(L)$ 
is a totally geodesic cross section that is orthogonal to the rulings.

\begin{theorem}\po\label{main1}Let $g\colon L^2\to\R^{n+2}$, $n\geq 4$, 
be a $2$-isotropic substantial surface. Then the associated  immersion 
$F_g\colon M^n\to\R^{n+2}$ is an austere $(n-2)$-ruled submanifold with 
complete rulings that has rank $\rho=4$ on an open dense subset of $M^n$. 
Moreover, the surface $j\colon L^2\to M^n$ is the unique totally geodesic 
cross section that is orthogonal to the rulings. Furthermore, the metric 
of $M^n$ is complete if and only if $L^2$ is complete.

Conversely, let  $F\colon M^n\to\R^{n+2}$, $n\geq 4$, be an austere 
$(n-2)$-ruled isometric immersion that has rank $\rho=4$ on an open dense 
subset of $M^n$. If there exists a totally geodesic global 
cross section $j\colon L^2\to M^n$ orthogonal to the rulings, then
the surface $g=F\circ j\colon L^2\to\R^{n+2}$ is $2$-isotropic and 
$F$ can be parametrized by $F_g$.
\end{theorem}
 
Assume that $L^2$ is simply connected. By Theorem $4$ in \cite{DV} 
there is a one-parameter family of minimal isometric immersions 
$F_\theta$ for $\theta\in [0,\pi)$ with $F_0=F_g$ such that each 
$F_\theta$ is also austere  carrying the same rulings and relative 
nullity subspaces as $F_g$. Consequently, we have the isometric immersions 
in higher codimension
$$
G=(\cos\varphi F_0,\sin\varphi F_{\pi/2})\colon M^n\to\R^{n+2}\oplus\R^{n+2}
\equiv\R^{2n+4},\;\;\varphi\in [0,\pi],
$$
that are also austere with the same rulings and relative nullity subspaces.
\vspace{1ex}

The following result analyzes when the submanifold $M^n$ above is 
Kaehler, which turns out to be always the case for $n=4$. On the other 
hand, we see that the property of being Kaehler is exceptional for higher 
even dimensions.

\begin{theorem}\po\label{main2} Let $F_g\colon M^n\to\R^{n+2}$, 
$n\geq 4$, be the austere $(n-2)$-ruled submanifold associated to 
a $2$-isotropic  substantial surface $g\colon L^2\to\R^{n+2}$. Then  
$M^n$ is Kaehler if and only if $g$ is holomorphic. In addition $F_g$ 
in the Kaehler case is never holomorphic.
\end{theorem}

If $M^n$ above is Kaehler and simply-connected, being $F_g$ not
holomorphic it follows from a result in \cite{DG0} that $F_g$ admits an 
non-trivial associated one-parameter family of isometric minimal 
immersions. It can be shown that this family coincides with the one 
discussed after Theorem \ref{main1}.

\section{The proofs}

Let $g\colon L^2\to\R^{n+2}$, $n\geq 4$, be a substantial $1$-isotropic 
isometric immersion. Hence the surface is minimal and the first ellipse 
of curvature is a circle at all points. The minimality condition yields  
that the normal bundle of $g$ splits along an open dense subset of $L^2$ 
as the orthogonal sum 
\be\label{bundles}
N_gL=N_1^g\oplus N_2^g\oplus\dots\oplus N_m^g,\;\;\; m=[(n+1)/2],
\ee
of the higher normal bundles and these have rank two except possible the 
last one that has rank one if $n$ is odd. Given an orthonormal tangent
frame $\{e_1,e_2\}$ we have 
$$
N^g_k(p)=\spa\left\{\a_g^{k+1}(e_1,\ldots,e_1,e_1)(p),\;
\a_g^{k+1}(e_1,\ldots,e_1,e_2)(p)\right\}
$$
at $p\in L^2$.
Here $\a_g^2=\a_g\colon TL\times TL\to N_gL$ is the second fundamental 
form of $g$ and $\a_g^s\colon TL\times\cdots\times TL\to N_gL$, 
$s\geq 3$, is the $s^{th}$-fundamental form defined inductively by
$$
\a_g^s(X_1,\ldots,X_s)=\left(\nabla^\perp_{X_s}\ldots
\nabla^\perp_{X_3}\a_g(X_2,X_1)\right)^\perp
$$
where $(\;\;)^\perp$ denotes the projection onto the normal 
complement of $N^g_1\oplus\cdots\oplus N^g_{s-2}$.

The $k^{th}$-order ellipse of curvature 
${\cal E}_k^g(p)\subset N_k^g(p)$ at $p\in L^2$ is
$$
{\cal E}_k^g(p)=\{\a_g^{k+1}(e(\theta),\ldots,e(\theta))(p): e(\theta)
=\cos\theta e_1+\sin\theta e_2\;\;\mbox{and}\;\;\theta\in [0,2\pi)\}.
$$
Then ${\cal E}_k^g(p)$ is indeed an ellipse and is a circle if and 
only if the vectors $$\a_g^{k+1}(e_1,\ldots,e_1,e_1)(p),\; 
\a_g^{k+1}(e_1,\ldots,e_1,e_2)(p)$$ are orthogonal with equal norm.
\vspace{1ex}

\noindent\emph{Proof of Theorem \ref{main1}:}
The minimal submanifold $F_g\colon M^n\to\R^{n+2}$ parametrized by 
\eqref{param} is  $(n-2)$-ruled   of rank four. Its tangent bundle 
splits orthogonally as  
$$
TM=\mathcal{H}\oplus\mathcal{V}
$$
where ${\cal H}$ is the tangent distribution orthogonal to the rulings and
$\mathcal V=\ker\pi_*$ is the vertical bundle of the submersion $\pi$.
Then $j_*T_pL={\cal H}(j(p))$ at every point $p\in L^2$ and the fibers 
of $\mathcal V$ form the distribution tangent to the rulings. 
We also have the orthogonal splitting
$$
\mathcal{V}=\mathcal{V}^1\oplus\mathcal{V}^0
$$ 
where $\mathcal{V}^1$ is identified with the fibers of $N_2^g$ and
$\mathcal{V}^0$ is identified with $N_3^g\oplus\cdots\oplus N_m^g$ and  
are the relative nullity subspaces of $F_g$.

Let $e_1,\ldots,e_{n+2}$ be an orthonormal frame such that
$N_r^g=\spa\{e_{2r+1},e_{2r+2}\}$. Denote 
$$
\o^k_{ij}=\<\nab_{e_k}e_i,e_j\>,\;1\leq k\leq 2,\;1\leq i,j\leq n+2
$$
where $\nab$ stands for the connection in $\R^{n+2}$. By assumption 
$$
\a_{11}^2=\a_g^2(e_1,e_1)=\kappa_1 e_3\;\;\mbox{and}
\;\;\a_{12}^2=\a_g^2(e_1,e_2)=\kappa_1 e_4.
$$
Thus
\be
\begin{aligned}\label{tres}
\a_{111}^3=(\nap_{e_1}\a^2_{11})_{N_2^g}=\kappa_1(\nap_{e_1}e_3)_{N_2^g}
&=\kappa_1(\o_{35}^1e_5+\o_{36}^1e_6)=\kappa_1(a_1e_5+b_1e_6)\vspace{1ex}\\
\a_{112}^3=(\nap_{e_1}\a^2_{12})_{N_2^g}=\kappa_1(\nap_{e_1}e_4)_{N_2^g}
&=\kappa_1(\o_{45}^1e_5+\o_{46}^1e_6)=\kappa_1(a_2e_5+b_2e_6).
\end{aligned}
\ee

As shown in \cite{DV} there is an orthonormal tangent frame 
$E_i$, $1\leq i\leq n$, such that
\be\label{es}
\mathcal{H}=\spa\{E_1,E_2\},\;\;\mathcal{V}^1=\spa\{E_3,E_4\}
\;\;\mbox{and}\;\;\mathcal{V}^0=\spa\{E_5,\ldots,E_n\}
\ee
where $E_3,E_4$ are taken
constant in each ruling and $F_*E_j=e_{j+2}$, $3\leq j\leq n$.  
Then the submanifold can be parametrized as 
$$
F_g=g+\sum_{j=1}^{n-2}t_jE_{j+2}
$$
where $t_1,\ldots,t_{n-2}\in\R$. Moreover, there is an orthogonal 
normal frame $\xi,\eta$ satisfying $\|\xi\|=\Omega=\|\eta\|$ with 
$\Omega\in C^\infty(M)$ such that the shape operators of 
$F_g$ vanish on $\mathcal{V}_0$ and restricted to 
$\mathcal{H}\oplus\mathcal{V}^1$ have the form
\be\label{ssf}
A_{\xi}=\begin{bmatrix}
\kappa_1 + h_1&h_2&r_1&s_1\\
h_2&-\kappa_1-h_1&r_2&s_2\\
r_1&r_2&0&0\\
s_1&s_2&0&0&
\!\!\!\!\end{bmatrix},\;\;\;
A_{\eta}=\begin{bmatrix}
h_2&\kappa_1-h_1&r_2&s_2\\
\kappa_1-h_1&-h_2&-r_1&-s_1\\
r_2&-r_1&0&0\\
s_2&-s_1&0&0&
\!\!\!\!\end{bmatrix}.\\
\ee
Moreover, $r_j=-a_j/\Omega$, $s_j=-b_j/\Omega$, $j=1,2$,
with $\kappa_1,a_1,a_2,b_1,b_1\in C^\infty(L)$ 
whereas
$$
h_j=\frac{1}{\Omega^2}\left(t_1D_1^j+\cdots+t_4D_4^j\right),
\;\;j=1,2,
$$
where $D_i^j\in C^\infty(L)$, $1\leq i\leq 4$, and 
$t_1,\ldots,t_4\in\R$ are independent parameters.

We obtain from \eqref{ssf} that if $F_g$ is austere then the 
coefficients of the terms of third order of the characteristic 
polynomials of both shape operators have to vanish. From this
it turns out that austerity implies that
$$
2(r_1r_2+s_1s_2)h_2+(r_1^2+s_1^2-r_2^2-s_2^2)(\kappa_1+h_1)=0
$$
and
$$
2(r_1r_2+s_1s_2)(h_1-\kappa_1)-(r_1^2+s_1^2-r_2^2-s_2^2)h_2=0.
$$
It follows that austerity yields
\be\label{uno}
r_1r_2+s_1s_2=0\;\;\;\mbox{and}\;\;\;r_1^2+s_1^2=r_2^2+s_2^2
\ee
that is equivalent to
\be\label{dos}
a_1a_2+b_1b_2=0\;\;\;\mbox{and}\;\;\;a_1^2+b_1^2=a_2^2+b_2^2.
\ee
Using \eqref{tres} it follows that
\be\label{four}
\<\a_{111}^3,\a_{112}^3\>=\kappa_1^2(a_1a_2+b_1b_2)=0
\ee
and
\be\label{five}
\|\a_{111}^3\|^2=\kappa_1^2(a_1^2+b_1^2)
=\kappa_1^2(a_2^2+b_2^2)
=\|\a_{112}^3\|^2,
\ee
hence $g$ is $2$-isotropic.

To prove the converse, we have to verify that if \eqref{uno} 
holds then the coefficient of the term of third order of the 
characteristic polynomial of $A_{\cos\varphi\xi+\sin\varphi\eta}$ 
vanishes for any $\varphi\in [0,2\pi]$. In this case, since 
\eqref{four} and \eqref{five} hold we can choose $e_5$ and $e_6$ 
collinear with $\a_{111}^3$ and $\a_{112}^3$, respectively, and 
the remaining of the proof is just a long but straightforward computation.
\vspace{2ex}\qed

In the sequel, we will be dealing with the case when $M^n$ is a 
Kaehler manifold.
\vspace{2ex}

 Let $F=F_g\colon M^n\to\R^{n+2}$, $n=2m\geq 4$, be a minimal 
$(n-2)$-ruled submanifold associated to a $1$-isotropic oriented 
surface $g\colon L^2\to\R^{n+2}$. The orientation of $L^2$ induces 
an orientation on each plane vector bundle $N_k^g$ in \eqref{bundles} 
given by the ordered pair 
$$
\a_g^{k+1}(e_1,\ldots,e_1,e_1),\;\a_g^{k+1}(e_1,\ldots,e_1,e_2)
$$
where  $\{e_1,e_2\}$ is a positively oriented tangent frame. 
Then let the orthonormal frame $e_1,\ldots,e_{n+2}$ be such that
the pairs $e_{2r+1},\;e_{2r+2}$ spanning $N_r^g$ are positively
oriented. Now define $T\colon TM\to TM$ with respect to the  
orthonormal frame $E_1,\ldots,E_n$ as in \eqref{es} by 
$$
T|_{\mathcal{H}\oplus\mathcal{V}^1}
=\begin{bmatrix}
0&1&0&0\\
-1&0&0&0\\
0&0&0&-1\\
0&0&1&0&
\!\!\!\!\end{bmatrix}
$$
and $T|_{\mathcal{V}^0}=I$. Thus $T$ leaves invariant the distributions 
tangent to the rulings.

\begin{lemma}\po  The following facts are equivalent:
\begin{itemize}
\item[(i)]  $\a_F(TX,Y)=\a_F(X,TY)$ for all $X,Y\in TM$. 
\item[(ii)] $g$ is $2$-isotropic.
\end{itemize}
\end{lemma}

\proof  We have that $(i)$ is equivalent to 
$$
A_\xi\circ T=-T\circ A_\xi\;\;\mbox{and}\;\;A_\eta\circ T=-T\circ A_\eta.
$$
It is straightforward to verify that the above is equivalent to 
\be\label{equiv}
a_1=b_2,\;a_2=-b_1.
\ee
Thus \eqref{dos} holds and $g$ is $2$-isotropic.  Conversely, if
$g$ is $2$-isotropic then the pair of orthogonal vectors with the same
norm $\a_{111}^3,\a_{112}^3$ is positively oriented. Hence, we can take
$\a_{111}^3=\kappa e_5$ and $\a_{112}^3=\kappa e_6$ and \eqref{equiv} 
holds.
\vspace{1ex}\qed

\noindent\emph{Proof of Theorem \ref{main2}:}
Assume that $g$ is holomorphic. Then
$$
\a_g^{s+1}(e_1.\ldots,e_1)=\kappa_se_{2s+1}\;\;\mbox{and}\;\;
\a_g^{s+1}(e_1.\ldots,e_1,e_2)
=\kappa_se_{2s+2},\;\;1\leq s\leq n/2.
$$
Moreover,  from \cite{DV} the connection forms 
$\o_{\a,\beta}=\<\nap e_\a,e_\beta\>$ satisfy
\be\label{o1}
\omega_{2s-1,2s+1}=\omega_{2s,2s+2}=\tau_s\omega_1,
\ee
and
\be\label{o2}
\omega_{2s-1,2s+2}=-\omega_{2s,2s+1}=\tau_s\omega_2
\ee
where $\omega_1, \omega_2$ are dual to $e_1,e_2$, respectively, 
and $\tau_s=\kappa_s/\kappa_{s-1}$ with $\kappa_0=1$, $1\leq s\leq n/2$. 

Let $E_1,\ldots,E_n$ be an orthonormal frame as in \eqref{es}.  
We have to show that the almost complex structure $J$ defined as 
$J|_{\mathcal{H}\oplus\mathcal{V}^1}
=T|_{\mathcal{H}\oplus\mathcal{V}^1}$
and $JE_{2i+1}=E_{2i+2}$, $JE_{2i+2}=-E_{2i+1}$, $i\geq 2$,
is parallel. That is,
$$
\<\n_{E_k}E_i,E_j\>=\<\n_{E_k}JE_i,JE_j\>,
$$
or equivalently,
$$
\<\nab_{E_k}F_*E_i,F_*E_j\>=\<\nab_{E_k}F_*JE_i,F_*JE_j\>,
\;\;k=1,2\;\;\mbox{and}\;\;1\leq i,j\leq n.
$$
Since $g$ holomorphic, we have from \cite{DV} that
$$
F_*E_1=\frac{1}{\Omega}(g_*e_1-\tau_2(t_1e_3+t_2e_4)),\;\; 
F_*E_2=\frac{1}{\Omega}(g_*e_2-\tau_2(t_2e_3-t_1e_4)).
$$
We only argue for nontrivial cases: 
\vspace{1ex}

\noindent Let $i=1$ and $j=2s+1,\;s\geq 2$. Then
$$
\<\nab_{E_k}F_*E_1,F_*E_{2s+1}\>
=-\<\nab_{E_k}F_*E_2,F_*E_{2s+2}\>
$$
$$
\iff \<g_*e_1-\tau_2(t_1e_3+t_2e_4),\nab_{E_k}e_{2s+3}\>
=-\<g_*e_2-\tau_2(t_2e_3-t_1e_4),\nab_{E_k}e_{2s+4}\>.
$$
$$
\iff t_1\o_{3,2s+3}+t_2\o_{4,2s+3}
=-t_2\o_{3,2s+4}+t_1\o_{4,2s+4}.
$$
The last equality holds trivially for $s\geq 2$ and 
by \eqref{o1} and \eqref{o2} for $s=1$.
The proof for the cases $i=1,2$ and $j=2s+1,2s+2$, 
$s\geq 1$, is similar.\vspace{1ex}

\noindent Let $i=2s+1$ and $j=2r+1$ with $r\neq s$. Then
$$
\<\nab_{E_k}F_*E_{2s+1},F_*E_{2r+1}\>
=\<\nab_{E_k}F_*E_{2s+2},F_*E_{2r+2}\>
\iff \o_{2s+3,2r+3}=\o_{2s+4,2r+4}
$$
where the last equality either holds trivially or follows from \eqref{o1}.  
The proof for the remaining cases is similar.
\medskip

Now let us assume that $M^n$ is Kaehler. Being $F_g$ is austere we have
that $g$ is $2$-isotropic.  Being $F_g$ minimal we have 
$$
\a_F(JX,Y)=\a_F(X,JY)
$$
for any $X,Y\in TM$.  It follows easily that the three subspaces in 
the decomposition $TM=\mathcal{H}\oplus\mathcal{V}^1\oplus\mathcal{V}^0$ 
are $J$-invariant. Hence
$$
J|_{\mathcal{H}\oplus\mathcal{V}^1}=\begin{bmatrix}
0&1&0&0\\
-1&0&0&0\\
0&0&0&-\e\\
0&0&\e&0&
\!\!\!\!\end{bmatrix}
$$
where $\e=\pm 1$.

Since $g$ is $2$-isotropic we have from \cite{DV} that
$$
A_{\xi}=\begin{bmatrix}
\kappa_1 + h_1&h_2&r&0\\
h_2&-\kappa_1-h_1&0&r\\
r&0&0&0\\
0&r&0&0&
\!\!\!\end{bmatrix},\;\;\;
A_{\eta}=\begin{bmatrix}
h_2&\kappa_1-h_1&0&r\\
\kappa_1-h_1&-h_2&-r&0\\
0&-r&0&0\\
r&0&0&0&
\!\!\!\!\end{bmatrix}\\
$$
on $\mathcal{H}\oplus\mathcal{V}_1$. That 
$A_{\xi}J+JA_{\xi}=0$ and $A_{\eta}J+JA_{\eta}=0$ hold
is equivalent to $\e=1$.

We define an isometry $J^\perp\colon N_gL\to N_gL$ by
$$
J^\perp e_3=-e_4,\;\;J^\perp e_4=e_3
\;\;\mbox{and}\;\;J^\perp e_{j+2}=F_*JE_j,\;j\geq 3.
$$
Then $J^\perp$ is an almost complex structure since
$$
J^\perp e_5=F_*JE_3=F_*E_4=e_6,\;\;
J^\perp e_6=F_*JE_4=-F_*E_3=-e_5
$$
and 
$$
(J^\perp)^2e_{j+2}=F_*J^2E_j=-F_*E_j=-e_{j+2}.
$$

We claim that $J^\perp$ is parallel with respect to the 
normal connection of $g$. Since $J$ is parallel, we have 
$$
\<\n_XE_i,JE_j\>=-\<\n_XJE_i,E_j\>
$$
which is equivalent to
$$
\<\nab_XF_*E_i,F_*JE_j\>=-\<\nab_XF_*JE_i,F_*E_j\>
$$
for any $X\in TM$. 
\vspace{2ex}

\noindent If $i,j\geq 3$, we have
$$
\<\nap_Xe_{i+2},J^\perp e_{j+2}\>
=-\<\nap_XJ^\perp e_{i+2},e_{j+2}\>
$$
which gives
\be\label{a}
\left((\nap_XJ^\perp)e_{i+2}\right)_{N_1^{g\perp}}=0,\;\;i\geq 3.
\ee

\noindent If $i\geq 3$ we have
$$
\<(\nap_XJ^\perp)e_{i+2},e_3\>
=-\<J^\perp e_{i+2},\nap_Xe_3\>+\<e_{i+2},\nap_Xe_4\>.
$$
Since $\nap_Xe_3,\nap_Xe_4\in N_1^g\oplus N_2^g$, we obtain
\be\label{unos}
\<(\nap_XJ^\perp)e_{i+2},e_3\>=0,\;\;i\geq 5,
\ee
and
\be\label{doss}
\<(\nap_XJ^\perp)e_5,e_3\>=\o_{36}(X)+\o_{45}(X)=0.
\ee
Similarly, we obtain
\be\label{cuatro}
\<(\nap_XJ^\perp)e_6,e_3\>=0.
\ee
If follows from \eqref{unos}, \eqref{doss} and 
\eqref{cuatro} that 
\be\label{b}
\<(\nap_XJ^\perp)e_{i+2},e_3\>=0,\;\;i\geq 3.
\ee
and in the same way that
\be\label{c}
\<(\nap_XJ^\perp)e_{i+2},e_4\>=0,\;\;i\geq 3.
\ee
It follows from \eqref{a}, \eqref{b} and \eqref{c} that
$$
(\nap_XJ^\perp)e_{i+2}=0,\;\;i\geq 3.
$$
The same type of arguments yield
$$
(\nap_XJ^\perp)e_i=0,\;\;i=3,4,
$$
and this proves the claim.

We have 
$$
J^\perp\a(e_1,e_1)=-\kappa_1e_4=-\a_{12}=-\a(Je_1,e_1)
$$
and 
$$
J^\perp\a(e_1,e_2)=\kappa_1e_3=\a_{11}=-\a(Je_2,e_1).
$$
Hence 
$$
J^\perp\a(X,Y)=-\a(JX,Y).
$$
Let $\tilde{J}\colon g^*T\R^{n+2}\to g^*T\R^{n+2}$ be defined as
$$
\tilde{J}|_{g_*TL}
=g_*\circ J\;\;\mbox{and}\;\;\tilde{J}|_{N_gL}=-J^\perp.
$$
It is now straightforward to verify that $\nab\tilde{J}=0$, 
that is, that $\tilde{J}$ is a complex structure in $\R^{n+2}$ 
that satisfies
$
\tilde{J}\circ g_*=g_*\circ J,
$
hence $g$ is holomorphic.

For the last statement, observe that if we had that $F_g$ is 
holomorphic then we would have in \eqref{ssf} that 
$A_\eta=\pm J\circ A_\xi$, and it is easy to verify that 
this cannot be the case.\qed

\vspace{.5in} {\renewcommand{\baselinestretch}{1}
\hspace*{-20ex}\begin{tabbing} \indent\= IMPA -- 
Estrada Dona Castorina, 110
\indent\indent\= Univ. of Ioannina -- Math. Dept.\\
\> 22460-320 -- Rio de Janeiro -- Brazil  \>
45110 Ioannina -- Greece\\
\> E-mail: marcos@impa.br \> E-mail: tvlachos@uoi.gr
\end{tabbing}}
\end{document}